\documentclass[12pt,reqno]{amsart}
\usepackage{color}
\usepackage[english]{babel}
\usepackage{amsfonts}
\usepackage{amsmath}
\usepackage{amsthm}
\usepackage{amssymb}
\usepackage[misc]{ifsym}
\usepackage{bbm}
\usepackage{mathrsfs}
\usepackage{esint}
\usepackage{faktor}
\usepackage[utf8]{inputenc}
\usepackage{afterpage}
\usepackage[left=2.8cm,right=2.8cm,top=2.8cm,bottom=2.8cm]{geometry}
\usepackage{graphicx}
\newcommand{\N}{\mathbb{N}}

\newcommand{\R}{\mathbb{R}}

\newcommand{\Div}{\mathrm{div} \, }
\renewcommand{\epsilon}{\varepsilon}
\renewcommand{\phi}{\varphi}

\newcommand{\dx}{\, {\rm d}x}

\newtheorem{lemma}{Lemma}[section]
\newtheorem{thm}[lemma]{Theorem}
\newtheorem{prop}[lemma]{Proposition}

\theoremstyle{definition}

\newtheorem{rmk}[lemma]{Remark}

\numberwithin{equation}{section}

\DeclareMathOperator*{\essinf}{ess \, inf}

\begin{document}
	
	\title[Continuous spectrum for a double-phase eigenvalue problem]{Continuous spectrum for a double-phase \\ unbalanced growth eigenvalue problem}
	
	\author[L. Gambera]{Laura Gambera}
	\address[L. Gambera]{Dipartimento di Matematica e Informatica, Universit\`a degli Studi di Catania, Viale A. Doria 6, 95125 Catania, Italy}
	\email{laura.gambera@unipa.it}
	\author[U. Guarnotta]{Umberto Guarnotta}
	\address[U. Guarnotta]{Dipartimento di Matematica e Informatica, Universit\`a degli Studi di Catania, Viale A. Doria 6, 95125 Catania, Italy}
	\email{umberto.guarnotta@unipa.it}
	\author[N.S. Papageorgiou]{Nikolaos S. Papageorgiou}
	\address[N.S. Papageorgiou]{Department of Mathematics, National Technical University of Athens, Zografou Campus, 15780 Athens, Greece}
	\email{npapg@math.ntua.gr}
	
	\maketitle
	
	\begin{abstract}
	We consider an eigenvalue problem for a double-phase differential operator with unbalanced growth. Using the Nehari method, we show that the problem has a continuous spectrum determined by the minimal eigenvalue of the weighted $p$-Laplacian.
	\end{abstract}
	
	\let\thefootnote\relax
	\footnote{\Letter \quad Umberto Guarnotta (umberto.guarnotta@unipa.it).}
	
	
	\section{Introduction and main result}
	
Let $\Omega\subseteq \R^N$ be a bounded domain with Lipschitz boundary $\partial \Omega.$
In this paper we study the following double-phase eigenvalue problem	\begin{equation}
		\label{prob}
		\tag{{${\rm P_{\lambda}}$}}
		\left\{
		\begin{alignedat}{2}
			-\Delta_p^a u-\Delta_q u &= \lambda a(x)|u|^{p-2}u \quad &&\mbox{in} \;\; \Omega, \\
			u &= 0 \quad &&\mbox{on} \;\; \partial\Omega, \\
		\end{alignedat}
		\right.
	\end{equation}
	where $1<q<p<N $ and $\lambda>0.$  Given $a\in L^{\infty}(\Omega)_{+}\setminus \{0\},$ by $\Delta_p^a$ we denote the weighted $p$-Laplace differential operator defined by
	$$\Delta_p^a u:= \Div (a(x)|\nabla u|^{p-2}\nabla u).$$
If $a\equiv 1$, then we recover the standard $p$-Laplacian. The differential operator in \eqref{prob} is not homogeneous and is related to the so-called \textit{double phase integral functional} defined by
	\begin{equation}
		J(u):=\int_{\Omega} (a(x)|u|^p+ |u|^q)\, \dx.
	\end{equation}
Let $\Theta(x,t)$ be the integrand in this integral functional, that is,
	\begin{equation}
	\Theta(x,t):=a(x)t^p+t^q \quad \text{for all}\,\, (x,t)\in\Omega\times[0,+\infty).
	\end{equation}
	We do not assume that the weight function $a(\cdot)$ is bounded below by a positive constant, that is, we do not require that $\essinf_{\Omega}a>0$. Hence the integrand $\Theta (x, t)$ exhibits unbalanced growth, namely,
	\begin{equation}
	t^q\le\Theta(x,t)\le c_0(t^p+1) \quad \text{for all}\,\, (x,t)\in\Omega\times[0,+\infty),
\end{equation} 
for a suitable $c_0>0$.

Such integral functionals were first investigated by Marcellini \cite{PM1,PM2} and Zhikov \cite{Z} in the context of problems of calculus of variations and nonlinear elasticity theory. In the last decade the interest for such problems was revived and the regularity properties of minimizers of such functionals was investigated. We mention the important works of Baroni-Colombo-Mingione \cite{BCM}, Marcellini \cite{PM3,PM4} and Ragusa-Tachikawa \cite{RT}. We also recall the survey papers of Mingione-R\u{a}dulescu \cite{MR}, Papageorgiou \cite{P}, and  R\u{a}dulescu \cite{R}. 

A global regularity theory, similar to that for balanced growth problems (vide Lieberman \cite{L}), remains elusive and  this is a serious handicap in the study of unbalanced double phase problems. In the recent years existence and multiplcity results were proved for various types of unbalanced double phase equations. We mention the works of Colasuonno-Squassina \cite{CS}, Deregowska-Gasi\'nski-Papageorgiou \cite{DGP}, Gasi\'nski-Papageorgiou \cite{GP}, Gasi\'nski-Winkert \cite{GW}, Liu-Dai \cite{LD}, Liu-Papageorgiou \cite{LP}, Papageorgiou-Pudelko-R\u{a}dulescu \cite{PPR}, Papageorgiou-R\u{a}dulescu-Zhang \cite{PRZ}, Papageorgiou-Vetro-Vetro \cite{PVV}.

In this paper we show that problem \eqref{prob} has a continuous spectrum. More precisely, we prove that, for all $\lambda >\bar{\lambda}$, \eqref{prob} has a nontrivial solution (i.e., an eigenfunction), while for $\lambda\in [0, \bar{\lambda}]$  there is no solution. The critical value $\bar{\lambda}>0$ is identified in terms of the spectrum of $\Delta_p^a$ as derived in the recent work of Papageorgiou-Pudelko-R\u{a}dulescu \cite{PPR}. Similar results for the balanced double phase operator (that is, the $(p,q)$-Laplacian) can be found in the recent works of Bhattacharaya-Emamizadeh-Farjudian \cite{BEF}, Papageorgiou-R\u{a}dulescu \cite{PR}, and Papageorgiou-Vetro-Vetro \cite{PVV2}.\\
We point out that double phase equations provide models describing strongly anisotropic materials. The modularing coefficient $a(x)\ge0$ dictates the geometry of the composite made by two different materials, one with hardening exponent $p$ and the other with hardening exponent $q$.

\section{Preliminaries}
The unbalanced growth of $\Theta(x, \cdot)$ leads to a different functional framework, which is based not on the classical Sobolev spaces but on Musielak–Sobolev-Orlicz spaces. A comprehensive presentation about the theory of these spaces can be found in the book of Harjulehto-H\"{a}sto \cite{HH}. The lack of a global (i.e., up to the boundary) regularity theory for unbalanced double phase problems prevents the use of many powerful tools usually exploited in balanced double phase problems. For this reason our approach is  based on the Nehari method.

Given any domain $ \Omega \subseteq \R^N $, the space $ C^\infty_c(\Omega) $ consists of the compactly supported test functions in $ \Omega $. Denoted by $ L^{0}(\Omega) $ the set of all measurable functions $u:\Omega \to \R$, then the Musielak–Orlicz spaces are defined as
	\begin{equation*}
		L^{\Theta}(\Omega) := \{u \in L^{0}(\Omega): \,\rho_{\Theta}(u)< \infty\},
	\end{equation*}
where $$\rho_{\Theta}(u):=\int_{\Omega}\Theta(x,|u|) \dx=\int_{\Omega}(a(x)|u|^p+|u|^q)\dx.$$	We equip $L^{\Theta}(\Omega)$ with the so-called \textit{Luxemburg norm}, defined by $$ ||u||_{\Theta}:= \inf \biggl\{ \lambda>0: \rho_{\Theta}\biggl(\frac{u}{\lambda }\biggr)\le 1\biggr\}.$$ 
Endowed with this norm, $L^{\Theta}(\Omega)$ becomes a separable, reflexive Banach space. The Musielak–Sobolev-Orlicz space $W^{1,\Theta}
(\Omega)$ is defined by 
$$W^{1,\Theta}(\Omega):=\{ u \in L^{\Theta}(\Omega): |\nabla u|\in L^{\Theta}(\Omega)\}.$$
This space is equipped with the norm 
$$||u||_{1,\Theta}:= ||u||_{\Theta}+ ||\nabla u||_{\Theta}.$$
The space $ W^{1,\Theta}_0(\Omega) $ contains exactly the functions $ u \in W^{1,\Theta}(\Omega) $ such that $ u = 0 $ on $ \partial \Omega $ in the sense of traces; if $ \Omega $ is bounded, we endow $ W^{1,\Theta}_0(\Omega) $ with the standard equivalent norm given by the Poincaré inequality, that is,
\begin{equation}
	\label{sobolevnorm}
	\|u\|_{W^{1,\Theta}_0(\Omega)} := \|\nabla u\|_{L^\Theta(\Omega)}.
\end{equation}
The spaces $W^{1,\Theta}(\Omega)$ and $W^{1,\Theta}_0(\Omega)$ are separable, reflexive, Banach spaces (vide \cite[Proposition 2.14]{CS}). Let $\overline{\Theta}$ be the Young conjugate of $\Theta$, that is,
$$\overline{\Theta}(t):=\max_{s\ge0}\{st-\Theta(s)\}\quad \forall t \ge 0.$$ We indicate with $W^{-1, \overline{\Theta}}(\Omega)$ the topological dual of $W^{1,\Theta}_0(\Omega)$. Given $u\in W^{1,\Theta}_0(\Omega)$, by $\rho_\Theta(\nabla u)$ we mean $\rho_\Theta(|\nabla u|)$, and similarly for $\Theta_0$. \\
Finally, we denote by $C^{0,1}(\overline{\Omega})$  the space of Lipschitz functions $u:\overline{\Omega}\to \R $ and by $A_{p}$ the $p$-Muckenhoupt class of weight functions (see Cruz-Uribe and Fiorenza \cite[p.152]{CF}, as well as Harjulehto-H\"{a}sto \cite[p.106]{HH}).

We require that the differential operator $u\mapsto\Div (a(x)|\nabla u|^{p-2}\nabla u)+\Div (|\nabla u|^{q-2}\nabla u)$, where  $1<q<p<N$ and $\frac{p}{q}<1+\frac{1}{N}$, satisfies the following hypotheses:
\begin{itemize}
	\item[$ {\rm (H_1)} $]$a \in C^{0,1}(\overline{\Omega}) \cap A_{p},$
	\item[$ {\rm (H_2)} $] $a(x)>0$ for all $x \in \Omega.$
\end{itemize}
We denote by ${\rm (H)}$ the set of assumptions ${\rm (H_1)}$-${\rm (H_2)}$.
\begin{rmk}
The condition on the exponents $p,q$ implies that they cannot be far apart and also leads to $p<q^*=\frac{Nq}{N-q}$, which in turn guarantees compact embeddings for some relevant spaces. Moreover, under this condition on the exponents, the Poincaré inequality is available (vide \cite[Proposition 2.18]{CS}).
\end{rmk}

		Now we present some tools that will be used in the sequel. The following results
		(cf. \cite[Chapter 6]{HH}) furnish some useful embeddings and highlight some relation between the norm $||\cdot||_{\Theta}$ and the modular function $\rho_{\Theta}(\cdot)$. For their proofs, we refer to \cite[Proposition 2.15 and Proposition 2.17]{CS} and \cite[Proposition 2.1]{LD}, respectively.
\begin{prop}
	\label{embeddingresult}
	If hypothesis ${\rm (H)}$ holds, then
	\begin{itemize}
		\item[$ {\rm (a)} $] $L^{\Theta}(\Omega)\hookrightarrow L^{r}(\Omega )$ and $W^{1,\Theta}_{0}(\Omega)\hookrightarrow W^{1,r}_{0}(\Omega )$ continuously for all $r\in[1,q];$
		\item[$ {\rm (b)} $] $W^{1,\Theta}_{0}(\Omega)\hookrightarrow  L^{r}(\Omega )$ continuously for all $r\in[1,q^*]$ and compactly for all $r\in [1, q^*)$;
		\item[$ {\rm (c)} $] $L^{p}(\Omega)\hookrightarrow L^{\Theta}(\Omega )$ continuously.
	\end{itemize} 
\end{prop}

\begin{prop}
	\label{relationrhonorma}
	Suppose ${\rm (H)}$. Then one has
	\begin{itemize}
		\item[$ {\rm (a)} $] $||u||_{L^\Theta(\Omega)}=\lambda \Leftrightarrow \rho_{\Theta}\big(\frac{u}{\lambda}\big)=1;$
		\item[$ {\rm (b)} $] $||u||_{L^\Theta(\Omega)}<1$ (resp. $=1,>1$) $\Leftrightarrow \rho_{\Theta}(u)<1$ (resp. $=1,>1$);
		\item[$ {\rm (c)} $] $||u||_{L^\Theta(\Omega)}<1$  $\Rightarrow ||u||_{L^\Theta(\Omega)}^{p}\le \rho_{\Theta}(u) \le ||u||_{L^\Theta(\Omega)}^{q};$
		\item[$ {\rm (d)} $] $||u||_{L^\Theta(\Omega)}>1$  $\Rightarrow ||u||_{L^\Theta(\Omega)}^{q}\le \rho_{\Theta}(u) \le ||u||_{L^\Theta(\Omega)}^{p};$	
		\item[$ {\rm (e)} $] $||u||_{L^\Theta(\Omega)}\to 0$ (resp. $\to \infty$) $ \Leftrightarrow \rho_{\Theta}(u) \to 0$ (resp. $\to \infty$).	
	\end{itemize}
\end{prop}
\begin{flushleft}
	Now, we consider $A_p^{a}, A_q: W^{1, \Theta}_0(\Omega)\to W^{-1, \overline {\Theta}}(\Omega)$ defined by
\end{flushleft}
\begin{equation*}
	\begin{split}
		\langle A_p^{a}(u), h\rangle &=\int_{\Omega}a(x)|\nabla u|^{p-2}\langle \nabla u, \nabla h\rangle  \dx ,\\\langle A_q(u), h\rangle &=\int_{\Omega}|\nabla u|^{q-2}\langle \nabla u, \nabla h\rangle  \dx,
	\end{split}
\end{equation*}

\begin{flushleft}
	and we set $V:= A_p^{a}+A_q:W^{1, \Theta}_0(\Omega)\to W^{-1, \overline {\Theta}}(\Omega).$
	This operator is bounded, continuous, strictly monotone, and of type ($S_{+}$) (vide \cite[Proposition 3.1]{LD}).
	Finally we consider also the integrand 	\begin{equation}
		\Theta_{0}(x,t)=a(x)t^p \quad \text{for all}\,\, (x,t) \in \Omega\times[0,+\infty).
	\end{equation}
\end{flushleft}
Like before, we can define $L^{\Theta_{0}}(\Omega)$ and $W^{1, \Theta_{0}}(\Omega)$ and we observe that they are separable, reflexive Banach spaces (see \cite[pp.52-66]{HH}). Moreover $W^{1, \Theta}_{0}(\Omega)\hookrightarrow W^{1, \Theta_{0}}_{0}(\Omega)$ continuously and densely. Finally from Papageorgiou-R\u{a}dulescu-Zhang \cite[Lemma 2]{PRZ}, we know that under hypothesis ${\rm (H)}$ we have 
\begin{equation}
	\label{compactembeddingsthetaO}
	W^{1,\Theta_{0}}_{0}(\Omega)\hookrightarrow  L^{\Theta_{0}}(\Omega ) \quad \text{compactly.}
\end{equation}

Now we recall some results regarding the spectrum of the following nonlinear eigenvalue problem:
\begin{equation}
	\label{eigenvalueprob}
	\left\{
	\begin{alignedat}{2}
		-\Delta_p^a u &= \lambda a(x)|u|^{p-2}u \quad &&\mbox{in} \;\; \Omega, \\
			u &= 0 \quad &&\mbox{on} \;\; \partial\Omega. \\
	\end{alignedat}
	\right.
\end{equation}
Using \eqref{compactembeddingsthetaO}, Papageorgiou, Pudelko, and R\u{a}dulescu \cite{PPR} showed that there exists a smallest eigenvalue $\hat{\lambda}_{1}>0$  for problem \eqref{eigenvalueprob}, which is isolated in the spectrum of \eqref{eigenvalueprob} and simple. Moreover, they provide the following variational characterization of $\hat{\lambda}_{1}>0$:
\begin{equation}
	\label{characterization eigenvalue hatlambda}
	\hat{\lambda}_{1}=\inf{\biggl\{\frac{\rho_{\Theta_{0}}(\nabla u)}{\rho_{\Theta_{0}} (u)}: u \in W^{1, \Theta_{0}}_{0}(\Omega), u\not=0 \biggr\}}.
\end{equation}
This characterization can be rewritten using the homogeneity as follows:
\begin{equation}
	\label{characterization eigenvalue hatlambda normalizzata}
	\hat{\lambda}_{1}=\inf{\biggl\{\rho_{\Theta_{0}}(\nabla u):  \rho_{\Theta_{0}} (u)}=1\biggr\}.
\end{equation}
Furthermore, \cite[Proposition 5]{PPR} ensures that the elements of the eigenspace associated with $\hat{\lambda}_{1}$ are bounded functions, that is, belong in $W^{1, \Theta_{0}}(\Omega)\cap L^{\infty}(\Omega)$ and have constant sign. In fact, $\hat{\lambda}_{1}$ is the only eigenvalue with eigenfunctions with costant sign. All the other eigenvalues have nodal (i.e., sign-changing) eigenfunctions.\\ 
In analogy with \eqref{characterization eigenvalue hatlambda}, we define the following number related to \eqref{prob}:
\begin{equation}
	\label{characteriz. firsteigenvaluelambda^*}
	\lambda^{*}=\inf{\biggl\{\frac{\rho_{\Theta_{0}}(\nabla u)+ \frac{p}{q}||\nabla u||_q^{q}}{\rho_{\Theta_{0}} (u)}: u \in W^{1, \Theta}_{0}(\Omega), u\not=0 \biggr\}.}
\end{equation}
\begin{prop}
	Under ${\rm (H)}$, one has that $\lambda^{*}=\hat{\lambda}_{1}>0$.
	\begin{proof}
		From \eqref{characterization eigenvalue hatlambda} and \eqref{characteriz. firsteigenvaluelambda^*} it is clear that $\hat{\lambda}_{1}\le \lambda^*$. Let $\hat{u}\in  W^{1, \Theta_{0}}_{0}(\Omega)$ be an eigenfunction related to the eigenvalue $\hat{\lambda}_{1}>0$. For any $t>0$, using 
		\eqref{characteriz. firsteigenvaluelambda^*} and \eqref{characterization eigenvalue hatlambda}, we have
		\begin{equation}
			\begin{split}
				\lambda^{*}&\le \frac{\rho_{\Theta_{0}}(\nabla (t \hat{u}))+ \frac{p}{q}||\nabla (t \hat{u})||_q^{q}}{\rho_{\Theta_{0}} (t \hat{u})} = \frac{\hat{\lambda}_{1}\rho_{\Theta_{0}}(t \hat{u})+ \frac{p}{q}||\nabla (t \hat{u})||_q^{q}}{\rho_{\Theta_{0}} (t \hat{u})}\\
				&= \hat{\lambda}_{1}+ \frac{p}{q} \frac{1}{t^{p-q}}||\nabla  \hat{u}||_q^{q}.\\
			\end{split}
		\end{equation} 
		Letting $t\to +\infty$ one deduces $\lambda^*\le \hat{\lambda}_{1} $, concluding the proof.
	\end{proof}
\end{prop}
\section{Continuous Spectrum}
For any $\lambda>0$, the energy functional associated with \eqref{prob} is  $\varphi_{\lambda}:W^{1,\Theta}_{0}(\Omega)\to \R$ defined as
$$\varphi_{\lambda}(u)=\frac{1}{p}\rho_{\Theta_{0}}(\nabla u)+\frac{1}{q}||\nabla u||_{q}^{q}-\frac{\lambda}{p}\rho_{\Theta_{0}}( u) \quad \forall u \in W^{1,\Theta}_{0}(\Omega).$$
Then $\varphi_{\lambda}\in C^{1}(W^{1,\Theta}_{0}(\Omega))$ and we have 
$$\langle\varphi_{\lambda}'(u),h\rangle=\langle V(u),h\rangle-\lambda \int_{\Omega} a(x)|u|^{p-2}uh \dx \quad \forall u,h \in W^{1,\Theta}_{0}(\Omega).$$
As we already mentioned, our approach is based on the Nehari method. So we introduce the Nehari manifold for $\varphi_{\lambda}$ defined by 
$$\mathcal{N}_{\lambda}=\{ u\in  W^{1,\Theta}_{0}(\Omega)\setminus\{0\}:\langle\varphi_{\lambda}'(u),u\rangle=0\}.$$
\begin{prop}
	Under ${\rm (H)}$, for all $\lambda> \hat{\lambda}_{1}$, we have that $\mathcal{N}_{\lambda}\not=\emptyset.$
	\begin{proof}
		According to \eqref{characterization eigenvalue hatlambda}, we can find $\tilde{u}\in W^{1, \Theta_{0}}_{0}(\Omega)\setminus \{0\}$ such that $$\rho_{\Theta_{0}}(\nabla \tilde{u})< \lambda \rho_{\Theta_{0}}(\tilde{u}).$$ 
		Moreover, since $ W^{1, \Theta}_{0}(\Omega)\hookrightarrow W^{1, \Theta_{0}}_{0}(\Omega)$ continuously and densely, besides the modular function is continuous, we can request,
		without loss of generality, that  $\tilde{u}\in W^{1, \Theta}_{0}(\Omega)\setminus \{0\}$. Recalling that  $ W^{1, \Theta}_{0}(\Omega)\hookrightarrow W^{1, q}(\Omega)$ continuously (see Proposition \ref{embeddingresult}), we set $$t_{0}:=\bigg[\frac{||\nabla \tilde{u}||_{q}^{q}}{\lambda \rho_{\Theta_{0}}(\tilde{u})-\rho_{\Theta_{0}}(\nabla \tilde{u})}\bigg]^{\frac{1}{p-q}}>0. $$
		Then we have {\begin{equation*}
				\begin{split}
					t_{0}^{p}[\lambda\rho_{\Theta_{0}}(\tilde{u})-\rho_{\Theta_{0}}(\nabla \tilde{u})]	=t_{0}^{q}||\nabla \tilde{u}||_{q}^{q}\quad &\Leftrightarrow \quad \lambda\rho_{\Theta_{0}}(t_{0}\tilde{u})-\rho_{\Theta_{0}}(\nabla (t_{0}\tilde{u}))=||\nabla(t_{0}\tilde{u})||_{q}^{q}\\ 
     &\Leftrightarrow\quad
					\langle \varphi_{\lambda}'(t_{0}\tilde{u}),t_{0}\tilde{u}\rangle =0.
				\end{split}
		\end{equation*}}
Hence $t_{0}\tilde{u}\in \mathcal{N}_{\lambda}$, as desired.
	\end{proof}
\end{prop}
\begin{flushleft}
	Now we want to show that  $\hat{\lambda}_1$ is a lower bound for the eigenvalues related to \eqref{prob}.
\end{flushleft}
\begin{prop}
\label{Proposion smallest eigenvalue}
Let ${\rm (H)}$ be satisfied. Then any $\lambda\in(0,\hat \lambda_1]$ is not an eigenvalue of $(P_\lambda)$.
\begin{proof}
Arguing by contradiction, suppose that $(\lambda,u)$ with $\lambda \in (0, \hat \lambda_1]$ is an eigenpair of $(P_\lambda)$, that is, \begin{equation*}
	\label {eigenvalue}
	\rho_{\Theta_{0}}(\nabla u) + \|\nabla u\|_q^q = \lambda \rho_{\Theta_0}(u).
\end{equation*} 
By \eqref{characterization eigenvalue hatlambda}, besides $u\neq 0$, we deduce
\begin{equation*}
\hat{\lambda}_1 \rho_{\Theta_0}(u) < \hat{\lambda}_1 \rho_{\Theta_0}(u) + \|\nabla u\|_q^q \leq\rho_{\Theta_{0}}(\nabla u) + \|\nabla u\|_q^q = \lambda \rho_{\Theta_0}(u),
\end{equation*}
whence
\begin{equation*}
(\hat \lambda_1 - \lambda) \rho_{\Theta_0}(u) \,<\, 0.
\end{equation*}
Since $\hat \lambda_1 \geq \lambda$ and $\rho_{\Theta_0}(u)>0$ we get a contradiction. Thus,  any $\lambda \in (0, \hat \lambda_1]$ is not an eigenvalue of $(P_\lambda)$.
\end{proof}
\end{prop}
\begin{prop}
	\label{coercivity}
Suppose ${\rm (H)}$. If $\lambda > \hat \lambda_1$, then $\phi_{\lambda}|_{\mathcal N_\lambda}$ is coercive.
\begin{proof}
Suppose on the contrary that there exists a sequence $\{u_n\} \subseteq \mathcal N_{\lambda}$ and a constant $M>0$ such that
\begin{equation}
\label{absurd condition}
\|u_n\|_{W^{1,\Theta}_0(\Omega)} \rightarrow \infty \;\, \text{and} \;\, \phi_{\lambda}(u_n) \leq M  \quad \forall  n \in \mathbb N.
\end{equation}
Fix any $ n \in \mathbb N$ and set  $y_n := \dfrac{u_n}{ \| u_n\|_{\Theta_0}}$. Since $u_n \in \mathcal N_\lambda$, we have
\begin{equation}
\label{(7)}
\rho_{\Theta_0}(\nabla u_n) + \| \nabla u_n\|_q^q = \lambda \rho_{\Theta_0}(u_n).
\end{equation}
Consequently,
\begin{equation}
	\label{yn}
\frac{\rho_{\Theta_0}(\nabla u_n)}{ \| u_n\|_{\Theta_0}^q} + \| \nabla y_n\|_q^q = \lambda \frac{\rho_{\Theta_0}(u_n)}{\| u_n\|_{\Theta_0}^q}.
\end{equation}
Using \eqref{absurd condition} and the definition of $\phi_\lambda$ we get

\begin{equation*}
\frac{1}{q} \|\nabla u_n \|_q^q \leq M + \frac{1}{p} \big[\lambda \rho_{\Theta_0}(u_n) - \rho_{\Theta_0}(\nabla u_n)\big].
\end{equation*}
Dividing by $\| u_n\|_{\Theta_0}^q$ and exploiting \eqref{yn}, one arrives at 
\begin{equation*}
\begin{aligned}
 \frac{1}{q} \| \nabla y_n\|_q^q &\leq \dfrac{M}{\| u_n\|_{\Theta_0}^q} + \frac{1}{p} \bigg[ \dfrac{\lambda \rho_{\Theta_0}(u_n)- \rho_{\Theta_0}(\nabla u_n)}{\| u_n\|^q_{\Theta_0}}\bigg]\\
  &= \dfrac{M}{\| u_n\|_{\Theta_0}^q} + \frac{1}{p} \|\nabla y_n \|_q^q.
\end{aligned}
\end{equation*}
Hence
\begin{equation}
\label{(8)}
 \bigg[\frac{1}{q} - \frac{1}{p}\bigg] \| \nabla y_n \|_q^q \leq \dfrac{M}{\| u_n\|_{\Theta_0}^q}.
\end{equation}
We claim that $y_n \rightarrow 0 \ \text{ in } W^{1,q}_0(\Omega)$. According to \eqref{(8)}, it is sufficient to prove that $\| u_n\|_{\Theta_0} \rightarrow \infty$. To this end, we firstly notice that \eqref{absurd condition} and Proposition \ref{relationrhonorma} entail $\rho_{\Theta}(\nabla u_n) \rightarrow \infty$. 
 Moreover, from \eqref{(7)} we infer 
\begin{equation}
	\label{relation GW}
	\rho_{\Theta}(\nabla u_n)\, \le\, \lambda \rho_{\Theta_0}(u_n)\ \text{ for all } n \in \mathbb N,
\end{equation}
 so we have $\rho_{\Theta_{0}}(u_n) \rightarrow \infty$ and consequently $\| u_n\|_{\Theta_0} \rightarrow \infty$. Accordingly, \begin{equation}
	\label{yntendsto 0in w}
	y_n \rightarrow 0 \ \text{ in } W^{1,q}_0(\Omega).
\end{equation} 
Taking in account Proposition \ref{relationrhonorma} and the definition of $y_n$, \eqref{yn} can be rewritten as  
\begin{equation}
\begin{split}
\| \nabla y_n\|^q_q &= \| u_n\|_{\Theta_0}^{p-q}\big[ \lambda \rho_{\Theta_0}(y_n) - \rho_{\Theta_0}(\nabla y_n)\big] \\
& = \| u_n\|_{\Theta_{0}}^{p-q} \big[ \lambda - \rho_{\Theta_0}(\nabla y_n)\big]\ \quad \forall n\in \N,
\end{split}
\end{equation}
producing
\begin{equation*}
\rho_{\Theta_0}(\nabla y_n) = \lambda - \frac{\| \nabla y_n\|^q_q}{\| u_n\|_{\Theta_{0}}^{p-q}}.
\end{equation*}
Thus, recalling \eqref{yntendsto 0in w} and $\| u_n\|_{\Theta_0} \rightarrow \infty$, the sequence $\{y_n\}$ is bounded in $ W^{1,\Theta_{0}}_0(\Omega)$.
Hence, applying \eqref{compactembeddingsthetaO} and \eqref{yntendsto 0in w}, we obtain
$$
y_n \rightarrow 0 \text{ in } L^{\Theta_0}(\Omega),
$$
which contradicts the fact that $\| y_n\|_{\Theta_0} = 1$ for all $n \in \mathbb N$. This proves that $\phi_{\lambda}|_{\mathcal N_{\lambda}}$ is coercive.
\end{proof}
\end{prop}
Define $m_\lambda:=\inf_{\mathcal N_\lambda} \phi_\lambda$. We want to prove that, for all $\lambda > \hat \lambda_1$, $m_\lambda$ is a minimum for $\varphi_\lambda.$
\begin{prop}
\label{mlambda>0}
Suppose that  ${\rm (H)}$ holds. If $\lambda > \hat \lambda_1$, then  $m_\lambda > 0$.
\begin{proof}
Let $u\in \mathcal{N}_\lambda$. Then 
\begin{equation}
	\label{(11)}
	\rho_{\Theta_0}(\nabla u) + \| \nabla u\|_q^q = \lambda \rho_{\Theta_0}(u).
\end{equation}
Since $p>q$ and $u\not=0$, using \eqref{(11)} we get
\begin{equation}\label{varphimin}
	\begin{split}
		\varphi_{\lambda}(u)&=\frac{1}{p}\rho_{\Theta_{0}}(\nabla u)+\frac{1}{q}||\nabla u||_{q}^{q}-\frac{\lambda}{p}\rho_{\Theta_{0}}( u)\\
		&=\bigg[\frac{1}{q}-\frac{1}{p}\bigg]||\nabla u||_{q}^{q}>0.
	\end{split}
\end{equation}
Thus $m_\lambda\ge0$. In order to prove the statement, we argue by contradiction and suppose that $m_\lambda=0$. Consider a minimizing sequence $\{u_n\}\subseteq\mathcal{N}_\lambda$, i.e.,
\begin{equation}
	\label{min=0}
	\varphi_\lambda (u_n)\downarrow m_\lambda=0 \; \text{as}\; n\to \infty.
\end{equation}
Notice that \eqref{varphimin} and \eqref{min=0} imply
\begin{equation}
	\label{minq}
	u_n \to 0 \quad \text{in}\; W^{1,q}_0(\Omega).
\end{equation} Moreover, Proposition \ref{coercivity} ensures that $\{u_n\}\subseteq W^{1,\Theta}_0(\Omega)$ is bounded. So, reasoning up to subsequences, by reflexivity of $ W^{1,\Theta}_0(\Omega)$, there exists $\hat{u} \in W^{1,\Theta}_0(\Omega)$ such that
\begin{equation}
	\label{minweaktheta}
	u_n  \rightharpoonup \hat{u} \quad \text{in}\; W^{1,\Theta}_0(\Omega).
\end{equation}
From \eqref{minq}, \eqref{minweaktheta}, and Proposition \ref{embeddingresult} we infer $\hat{u}=0$.
Hence, exploiting the compactness of the
embedding $W^{1,\Theta}_{0}(\Omega)\hookrightarrow  L^{\Theta}(\Omega )$, we get
\begin{equation}	
	\label{limtheta L}
	u_n  \rightharpoonup 0 \quad \text{in}\; W^{1,\Theta}_0(\Omega)\quad \text{and} \quad	u_n \to 0 \quad \text{in}\; L^{\Theta}(\Omega).
\end{equation}
Now for any $n\in \mathbb{N}$ we define $v_n:=\frac{u_n}{||u_n||_{\Theta}}.$
Then, repeating verbatim the arguments in the proof of Proposition \ref{coercivity} and using \eqref{characteriz. firsteigenvaluelambda^*}, we ensure that
\begin{equation}
	\begin{split}
		\| \nabla v_n\|^q_q &= \| u_n\|_{\Theta}^{p-q}\big[ \lambda \rho_{\Theta_0}(v_n) - \rho_{\Theta_0}(\nabla v_n)\big] \\
		&\le \| u_n\|_{\Theta}^{p-q} (\lambda -\hat{\lambda}_{1} ) \rho_{\Theta_0}(v_n)\\
		&\le  \| u_n\|_{\Theta}^{p-q} (\lambda -\hat{\lambda}_{1})\rho_{\Theta}(v_n)\\
        &= \| u_n\|_{\Theta}^{p-q} (\lambda -\hat{\lambda}_{1}).
	\end{split}
\end{equation} 
Since $u_n\to 0$ in $L^{\Theta}(\Omega)$, besides recalling that $p<q^*$, one has  
\begin{equation}\label{vn}
v_n \to 0  \quad \text{in} \; W^{1,q}_0(\Omega) \quad \text{and}\quad	v_n \to 0  \quad \text{in} \; L^{p}(\Omega).
\end{equation}
It follows that $$ 0\le \rho_{\Theta_{0}}(v_n) = \int_\Omega a(x)|v_n|^p \dx \le \left(\max_{\overline{\Omega}} a\right) \int_\Omega |v_n|^p \dx \to 0.$$
Thus, exploiting also \eqref{vn}, we obtain $\rho_{\Theta}(v_n)\to 0$. By Proposition \ref{relationrhonorma} we get $||v_n||_{\Theta}\to 0$, contradicting $||v_n||_{\Theta}=1$ for all $n\in \N$. This proves that $m_{\lambda}>0$.
\end{proof}
\end{prop}

\begin{prop}
	\label{eigenfunction}
Suppose that ${\rm (H)}$ is satisfied and $\lambda > \hat{\lambda}_1$. Then there exists $\hat{u} \in \mathcal{N}_{\lambda}$ such that $\varphi_{\lambda}(\hat{u}_{\lambda}) = m_{\lambda}$.
\begin{proof}
For any fixed $\lambda > \hat{\lambda}_1$, we consider $m_{\lambda}=\inf_{\mathcal N_\lambda}\varphi_{\lambda}.$
Reasoning as in Proposition \ref{mlambda>0}, we take a minimizing sequence $\{u_n\} \subseteq \mathcal N_{\lambda}$ and we notice that it is bounded in $W^{1,\Theta}_0(\Omega)$.
So, arguing as in Proposition \ref{mlambda>0}, there exists $\hat{u}_{\lambda} \text{ in } W^{1,\Theta}_0(\Omega)$ such that

\begin{equation}
\label{(17)}
u_n \rightharpoonup\hat{u}_{\lambda} \text{ in } W^{1,\Theta}_0(\Omega)\quad \text{and}\quad \ u_n \rightarrow \hat{u}_{\lambda}\text{ in } L^{\Theta}(\Omega).
\end{equation}
Since $u_n \in \mathcal N_{\lambda}$, we have
\begin{equation}
\label{(18)}
\rho_\Theta(\nabla u_n) = \rho_{\Theta_{0}}(\nabla u_n) + \|\nabla u_n\|_{q}^{q} = \lambda \rho_{\Theta_{0}}(u_n) \leq \lambda \rho_{\Theta}(u_n) \text{ for all } n \in \mathbb N.
\end{equation}
We prove $\hat{u}_{\lambda} \not=0$ arguing by contradiction. If $\hat{u}_{\lambda} =0$ then, passing to the limit in \eqref{(18)} via \eqref{(17)}, we get $\rho_\Theta(\nabla u_n)\to 0$, that is,
\begin{equation*}
	 u_n \rightarrow 0 \text{ in } W^{1, \Theta}_0(\Omega),
\end{equation*}
Since $W^{1,\Theta}_0(\Omega)\hookrightarrow W^{1,\Theta_0}_0(\Omega)\hookrightarrow L^{\Theta_0}(\Omega)$, by definition of $\phi_\lambda$ we have $\phi_\lambda(u_n)\to 0$, which entails $m_\lambda = 0$ by construction of $\{u_n\}$. This contradicts Proposition \ref{mlambda>0}, so $\hat{u}_{\lambda} \neq 0$. \\
Moreover, the weak lower sequential semi-continuity of $\varphi_{\lambda}$ implies that
$$
\varphi_{\lambda}(\hat{u}_{\lambda}) \leq \liminf\limits_{n \rightarrow \infty}\varphi_\lambda(u_n) = m_\lambda.
$$
From this inequality we get that if $\hat u_\lambda \in \mathcal N_\lambda$, then   $\varphi_\lambda(\hat u_\lambda) = m_\lambda$. Suppose by contradiction that $\hat u_\lambda \not \in \mathcal N_\lambda$ and consider the function $k_\lambda: [0,1]\to \R$ defined as
$$
k_\lambda(t) = \langle \varphi'_{\lambda}(t\hat u_\lambda),t\hat u_\lambda \rangle= \rho_{\Theta_0}(\nabla (t\hat u_\lambda)) + \|\nabla(t\hat u_\lambda)\|_q^q - \lambda \rho_{\Theta_0}(t\hat u_\lambda).
$$
Exploiting \eqref{(17)} and the weak lower sequential semi-continuity of both $\rho_{\Theta_0}(\cdot)$ and $\|\cdot\|_q$, besides recalling that $W^{1,\Theta}_0(\Omega) \hookrightarrow W^{1,\Theta_0}_0(\Omega)$, $W^{1,\Theta}(\Omega) \hookrightarrow L^q(\Omega)$, and $L^\Theta(\Omega) \hookrightarrow L^{\Theta_0}(\Omega)$, by $\hat u_\lambda \not \in \mathcal N_\lambda$ we have
\begin{equation*}
\rho_{\Theta_{0}}(\nabla \hat{u}_{\lambda}) + \|\nabla \hat{u}_{\lambda}\|_{q}^{q} < \lambda \rho_{\Theta_{0}}(\hat{u}_{\lambda})
\end{equation*}
that is,
\begin{equation}
\label{(21)}
k_\lambda (1) < 0.
\end{equation}
Moreover, according to \eqref{characterization eigenvalue hatlambda}, there exist  $c_1, c_2>0$ such that
\begin{equation}
k_\lambda(t) \geq t^q \|\nabla\hat u_\lambda \|_q^q - t^p(\lambda -\hat \lambda_1) \rho_{\Theta_0}(\hat u_\lambda) = c_1t^q -c_2 t^p.
\end{equation}
Since $q < p$, for $t \in (0,1)$ small we have 
\begin{equation}
\label{(22)}
k_\lambda(t) > 0.
\end{equation}
Hence we are in the position to apply Bolzano's Theorem, because of \eqref{(21)} and \eqref{(22)}. So we deduce that there exists $\hat t \in (0,1)$ such that
$$
\begin{aligned}
& k_\lambda(\hat t) =  \langle\varphi'_{\lambda}(\hat t\hat u_\lambda),\hat t\hat u_\lambda\rangle=0,
 \end{aligned}
$$
that is, $\hat t \hat u_\lambda \in \mathcal N_\lambda.$ Now we repeat the argument of \eqref{varphimin}  for $\hat t \hat u$. Since $0<\hat t<1$, exploiting \eqref{(17)} like above and recalling that $u_n\in\mathcal N_\lambda$, one has
\begin{equation}
\begin{aligned}
m_\lambda &\leq \varphi_\lambda(\hat t \hat u_\lambda) = \bigg[\frac{1}{q} - \frac{1}{p}\bigg]\hat t^q \| \nabla \hat u_\lambda \|_q^q \\
& < \bigg[\frac{1}{q}-\frac{1}{p}\bigg] \| \nabla \hat u_\lambda\|_q^q \\
& \leq \bigg[\frac{1}{q} - \frac{1}{p}\bigg] \liminf\limits_{n \rightarrow \infty} \| \nabla u_n\|_q^q \\
& = \liminf\limits_{n \rightarrow \infty}\varphi_\lambda(u_n) = m_\lambda.
\end{aligned}
\end{equation}
This is a contradiction. Therefore $\hat u_\lambda \in \mathcal N_\lambda$ and $m_\lambda = \varphi_\lambda(\hat u_\lambda)$.
\end{proof}
\end{prop}
Fix any $\lambda>\hat {\lambda}_1$ and let $K_{\varphi_{\lambda}}$ be the critical set of $\varphi_\lambda$, that is, 
$$
K_{\varphi_\lambda} = \{u \in W^{1,\Theta}_0(\Omega)\ :\ \varphi'_\lambda (u) = 0\}.
$$
Now we turn to show that $\hat u_\lambda \in K_{\varphi_\lambda}$, so that $\hat u_\lambda$ is an eigenfunction of problem $(P_\lambda)$. This implies that $(P_\lambda)$ has a continuous spectrum.
\begin{thm}
Let ${\rm (H)}$ be satisfied. Then every $\lambda > \hat \lambda_1$ is an eigenvalue of $(P_\lambda)$, with eigenfunction 
$$
\hat u_\lambda \in W^{1,\Theta}_0(\Omega)\cap L^{\infty}(\Omega), \quad \hat u_\lambda(x) > 0 \text{ for a.a. } x \in \Omega.
$$
\begin{proof}
Consider the function $\psi_{\lambda}: W^{1,\Theta}_0(\Omega) \rightarrow \mathbb R$ defined by
$$
\psi_{\lambda} (u) = \langle \phi'_\lambda(u),u\rangle = \rho_{\Theta_0}(\nabla u) + \| \nabla u\|_q^q - \lambda \rho_{\Theta_0}(u) \text{ for all } u \in W^{1,\Theta}_0(\Omega).
$$ 
Evidently $\psi_\lambda \in C^1(W^{1,\Theta}_0(\Omega))$. From Proposition \ref{eigenfunction}, we can find $\hat u_{\lambda}\not=0$ such that
$$
\varphi_\lambda(\hat u_\lambda) = \min\big\{\phi_\lambda(u) : u \in W^{1,\Theta}_0(\Omega), \psi_\lambda(u)= 0\big\}.
$$
Then, by the Lagrange multiplier rule \cite[p.442]{PRR}, there exists $\eta \in \mathbb R$ such that $$\phi'_\lambda(\hat u_\lambda) = \eta \psi'_\lambda(\hat u_\lambda).$$
To conclude, it is sufficient to show that $\eta=0$. Suppose, by contradiction, $\eta\not=0$. Recalling that $\hat u_\lambda \in \mathcal N_\lambda$, we obtain
\begin{equation*}
\eta \langle \psi'_\lambda(\hat u_\lambda), \hat u_\lambda \rangle = 0,
\end{equation*}
that is,
\begin{equation*}
\begin{aligned}
  \eta p\big[ \rho_{\Theta_0}(\nabla\hat u_\lambda)-\lambda \rho_{\Theta_0}(\hat u_\lambda)\big] +\eta q \| \nabla \hat u_\lambda\|^q_q &= 0.\end{aligned}
\end{equation*}
Since $q<p$ and $\eta\not=0$, we have that
\begin{equation*} 
 0=\eta p\big[ \rho_{\Theta_0}(\nabla\hat u_\lambda)-\lambda \rho_{\Theta_0}(\hat u_\lambda)\big] +\eta q \| \nabla \hat u_\lambda\|^q_q <\eta p\big[ \rho_{\Theta_0}(\nabla\hat u_\lambda) + \| \nabla\hat u_\lambda \|^q_q - \lambda \rho_{\Theta_0}(\hat u_\lambda) \big].
\end{equation*}
This is a contradiction since $\hat u_\lambda \in \mathcal N_\lambda$. Therefore $\eta = 0$, whence
\begin{equation*}
\begin{aligned}
\phi'_\lambda(\hat u_\lambda) = 0,
\end{aligned}
\end{equation*}
which implies $$\hat u_\lambda \in K_{\phi_\lambda}.$$
Hence $\hat u_\lambda$ is a nontrivial solution of problem $(P_\lambda)$. Moreover \cite[Theorem 3.1]{GW} guarantees that
$$
\hat u_\lambda \in W^{1,\Theta}_0(\Omega)\cap L^{\infty}(\Omega).
$$
Note that $\phi_\lambda(u) = \phi_\lambda(|u|)$ for all $u \in W^{1,\Theta}_0(\Omega)$.
Hence we can assume $\hat u_\lambda \geq 0$. Finally \cite[  Proposition 2.4] {PVV} implies that $\hat u_\lambda(x) > 0$ for a.a. $x \in \Omega$, which concludes the proof.
\end{proof}
\end{thm}	

\section*{Acknowledgments}

\noindent
U. Guarnotta is supported by: (i) PRIN 2017 `Nonlinear Differential Problems via Variational, Topological and Set-valued Methods' (Grant No. 2017AYM8XW) of MIUR; (ii) GNAMPA-INdAM Project CUP$\underline{\phantom{x}}$E55F22000270001; (iii) grant `PIACERI 20-22 Linea 3' of the University of Catania.

\noindent
\textbf{Conflict of interest statement.} The authors declare that they have no conflict of interest.

\end{document}